\input amstex\documentstyle{amsppt}  
\pagewidth{12.5cm}\pageheight{19cm}\magnification\magstep1
\topmatter
\title Superspecial representations of Weyl groups
\endtitle
\author G. Lusztig\endauthor
\address{Department of Mathematics, M.I.T., Cambridge, MA 02139}\endaddress
\thanks{Supported by NSF grant DMS-2153741}\endthanks
\endtopmatter   
\document

\define\Irr{\text{\rm Irr}}

\define\si{\sim}

\define\sqc{\sqcup}

\define\bin{\binom}
\define\op{\oplus}
   
\define\part{\partial}
\define\em{\emptyset}

\define\n{\notin}
\define\iy{\infty}
\define\m{\mapsto}
\define\do{\dots}

\define\sub{\subset}    
\define\bxt{\boxtimes}
\define\T{\times}
\define\ti{\tilde}
\define\nl{\newline}
\redefine\i{^{-1}}

\define\ot{\otimes}
\define\bbq{\bar{\QQ}_l}

\define\ind{\text{\rm ind}}

\define\sg{\text{\rm sgn}}

\define\a{\alpha}
\redefine\b{\beta}

\define\g{\gamma}

\define\r{\rho}
\define\s{\sigma}

\redefine\l{\lambda}

\define\x{\xi}

\redefine\G{\Gamma}

\define\CC{\bold C}

\define\NN{\bold N}

\define\QQ{\bold Q}
\define\RR{\bold R}

\define\ZZ{\bold Z}

\define\cf{\Cal F}

\define\sha{\sharp}

\head Introduction\endhead
\subhead 0.1\endsubhead
Let $W$ be a Weyl group. Let $\Irr(W)$ be the set of (isomorphism
classes of) irreducible representations over $\CC$ of $W$.

In this paper we define a subset $\Irr_{ssp}(W)$ of
the set of special representations of $W$. (See 1.3, 1.14.)

To do this we consider a connected reductive group $G$ with Weyl
group $W$ defined and split over a finite field $F_q$ with
group of rational points $G(F_q)$ and also an $F_q$-rational
structure on $G$ for which the Frobenius acts on the Weyl group
as opposition, with group of rational points $G(F_q)'$.
From \cite{L84} it is known that (for sufficiently large $q$)
there is a bijection $\r\m\r'$ from the set of
unipotent representationsof $G(F_q)$ to the set of
unipotent representations of $G(F_q)'$ such that $\pm1$ times
the dimension of $\r'$ (as a polynomial in $q$)
is obtained by replacing $q$ by $-q$ in the polynomial in $q$ which
gives the dimension of $\r$.
Our observation is that there is at most one $\r=\r_E$ in the
unipotent principal of $G(F_q)$ series which corresponds to a special
representation $E$ of $W$ and is such that $\r'$ is cuspidal for
$G'(F_q)$.
The $E$ obtained in this way are called superspecial.
The irreducible Weyl groups for which such $E$ exist are listed in
1.17.
It turns out that if $E$ is as above then $\dim(\r_E)$ is of the
form $q^{a_E}\sha(G(F_q))^*/(\pm c_EP_E(-q))$
where $a_E,c_E$ are independent of $q$,
$\sha(G(F_q))^*$  is the part prime to $q$ in $\sha(G(F_q))$  and
$P_E$ is a polynomial in $q$ with coefficients in $\NN$;
moreover, $P_E$ can be factored as a product of remarkably simple
polynomials (each one with coefficients in $\NN$).

Although the arguments
above were our motivation, our definition of superspecial
representations is actually
purely in terms of generic degrees and does not make use of the
groups $G(F_q)$ or $G'(F_q)$.

In \S2 we associate to a superspecial representation of $W$
(assumed to be irreducible) a constructible representation $Z_W$
of $W$ (or, in the nonsimply laced case, two such representations,
$Z_W,Z'_W$.) We will show elsewhere that using
these representations one can reconstruct
(without using algebraic geometry) 
the finite groups (products of symmetric groups) associated in
\cite{L84} to any special representation, which were used in
\cite{L84} to classify unipotent representations of reductive
groups over $F_q$.

In \S3 the definition of superspecial representations is extended
to finite noncrystallographic Coxeter groups.

In \S4 we associate (using \cite{L02}) to a superspecial
representation of $W$ (assumed to be irreducible with trivial
opposition) an elliptic conjugacy class of $W$, which we call
the superspecial conjugacy class.

\subhead 0.2. Notation\endsubhead
For $n\in\ZZ$ we set $[n/2]=n/2$ if $n\in2\ZZ$ and
$[n/2]=(n-1)/2$ if $n\in2\ZZ+1$.

\head 1. Definition of superspecial representations\endhead
\subhead 1.1\endsubhead
The set of simple reflections of $W$ is denoted by $I$. Let $r$ be
the number of orbits of the opposition involution $op:I@>>>I$.

Let $u$ be an indeterminate.
For $E\in\Irr(W)$ the generic degree $D_E(u)$ can be defined
in terms of the Iwahori-Hecke algebra of $W$ (see for example
\cite{AL,p.202}). (A priori, $D_E(u)$ is in the quotient field of
$\CC[u]$; in fact, it is in $\QQ[u]$.) It is known that    
$$D_E(u)=u^{a_E}\prod_{i=1}^r(u^{e_i+1}-1)/((-1)^{deg P_E}c_EP_E(-u))
\tag a$$
where $e_1,e_2,\do,e_r$ are the exponents of $W$, $a_E\in\NN$,
$c_E\in\{1,2,3,\do\}$ and $P_E(u)$ is a product of cyclotomic
polynomials.

\subhead 1.2\endsubhead
We say that $E$ is a {\it special representation} if $E$ appears in the
$a_E$-th symmetric power of the reflection representation of $W$.
(It then appears there with multiplicity $1$.)

Let $\Irr_{sp}(W)$ be the subset of $\Irr(W)$ consisting of special
representations. (This subset has been introduced in
\cite{L79a}.)

In \cite{L79}, the set $\Irr(W)$ has been partitioned into
{\it families};
it is known that any family contains a unique special representation.

For $E\in\Irr(W)$ we denote by $\g_E$ the largest integer such that
$(1+u)^{\g_E}$ divides the polynomial $D_E(u)$.

\proclaim{Theorem 1.3} (i) For any $E\in\Irr(W)$ we have $\g_E\le r$.

(ii) The set $\Irr_{ssp}(W):=\{E\in\Irr_{sp}(W);\g_E=r\}$ consists
of at most one element.

(iii) If $E\in\Irr_{ssp}(W)$ then $P_E(u)\in\NN[u]$. It has degree
$2a_E+\sha(I)$ and is a product of polynomials
of the form $1+u^s+u^{2s}+\do+u^{(l-1)s}$ where $s\in\{1,2,3,\do\}$
and $l$ is a prime number dividing $2c_E$. 
\endproclaim
To prove the theorem we can assume that $W$ is irreducible.
The various cases will be considered in 1.6-1.13.

We say that $E\in\Irr(W)$ is {\it superspecial} if it is contained in
$\Irr_{ssp}(W)$. We say that $W$ is {\it superspecial} if
$\Irr_{ssp}(W)\ne\em$; we then denote by $E_W$ the unique object of
$\Irr_{ssp}(W)$ and by $\cf_W$ the family of $\Irr(W)$ that contains
$E_W$.

\subhead 1.4\endsubhead
Let $X\sub\ZZ_{>0}$, $\sha(X)=m<\iy$. We set
$$n=\sum_{\l\in X}\l-\bin{m}{2}.$$ Note that $n>0$. Let
$X^0=X\cap2\NN,X^1=X\cap(2\NN+1)$, $m^0=\sha(X^0),m^1=\sha(X^1)$. Let
$$\a=\bin{m^0}{2}-\sum_{\l\in X^0}\l/2
+\bin{m^1}{2}-\sum_{\l\in X^1}\l/2+m^1/2+[n/2].\tag a$$
We show:

(b) We have $\a\le[n/2]$ with equality if and only if
$X=\{1,3,\do,2m-1\}$.
\nl
We have   
$$\bin{m^0}{2}-\sum_{\l\in X^0}\l/2\le
(1+2+\do+(m^0-1))-(2+4+\do+(2m^0))/2\le-m^0\le0$$
with the last $\le$ being $=$ if and only if $m^0=0$. We have
$$\bin{m^1}{2}+m^1/2-\sum_{\l\in X^1}\l/2\le
(1+3+\do+(2m^1-1))/2-(1+3+\do+(2m^1-1))/2=0$$
with $\le$ being $=$ if and only if $X^1=\{1,3,\do,2m^1-1\}$.
Since $m=m^0+m^1$, this proves (b).

\subhead 1.5\endsubhead
Let $(X,Y)\sub\NN\T\NN$ be such that $0\n Y$,
$x=\sha(X)<\iy$, $y=\sha(Y)<\iy$. We set $m=x+y$,
$$n=\sum_{\l\in X}\l+\sum_{\l\in Y}\l-[(m-1)^2/2]/2.$$
Let $$X^0=X\cap2\NN,X^1=X\cap(2\NN+1), x^0=\sha(X^0),x^1=\sha(X^1),$$
$$Y^0=Y\cap2\NN,Y^1=Y\cap(2\NN+1), y^0=\sha(Y^0),y^1=\sha(Y^1).$$ Let
$$\a=\bin{x^0}{2}+\bin{x^1}{2}+\bin{y^0}{2}+\bin{y^1}{2}+x^0y^1+x^1y^0
-\sum_{\l\in X}\l-\sum_{\l\in Y}\l+2[n/2].\tag a$$
We show:

(b) We have $\a\le2[n/2]$ with equality if and only if
$$(X,Y)=(\{0,2,4,\do,m-2\},\{1,3,\do,m-1\})$$ (with $m$ even)
or $$(X,Y)=(\{0,2,4,\do,m-1\},\{1,3,\do,m-2\})$$ (with $m$ odd).
\nl
We have $\a=\a^{01}+\a^{10}+2[n/2]$ where
$$\a^{01}=\bin{x^0+y^1}{2}-\sum_{\l\in X^0\sqc Y^1}\l,$$
$$\a^{10}=\bin{x^1+y^0}{2}-\sum_{\l\in X^1\sqc Y^0}\l.$$
We have
$$\a^{01}\le z$$
where 
$$z=(1+2+3+\do+(x^0+y^1-1))-(0+2+4+\do+(2x^0-2))-(1+3+\do+(2y^1-1).$$
If $X^0=\{0,2,4,\do,2x^0-2\},Y^1=\{1,3,\do,2y^1-1\}$, we have
$\a^{01}=z$; if this condition is not satisfied then $\a^{01}<z$
Now $z=0$ if $x^0=y^1$ or $x^0=y^1+1$ and $z<0$ in all other cases. 
We see that $\a^{01}=0$ if 
$$X^0=\{0,2,4,\do,2x^0-2\},Y^1=\{1,3,\do,2x^0-1\}$$ (for some $x_0\ge0$)
or if $$X^0=\{0,2,4,\do,2x^0-2\},Y^1=\{1,3,\do,2x^0-3\}$$
(for some $x_0\ge1$) and $\a^{01}<0$ in all other cases.

Similarly, $\a^{10}=0$ if 
$$Y^0=\{0,2,4,\do,2y^0-2\},X^1=\{1,3,\do,2y^0-1\}$$
(and $y_0\ge0$ must be $0$ since $0\n Y^0$)
or if $$Y^0=\{0,2,4,\do,2y^0-2\},X^1=\{1,3,\do,2y^0-3\}$$
(for some $y_0\ge1$, but this would imply $0\in Y^0$ which is
not the case); we have $\a^{10}<0$ in all other cases.
In other words, we have $\a^{10}=0$ if $Y^0=\em,X^1=\em$
(that is $y^0=x^1=0$) and $\a^{10}<0$ in all other cases.
When $x^0=y^1$, $y^0=x^1=0$ we have $x^0=y^1=m/2$ so that $m\in2\NN$.
When $x^0=y^1+1$, $y^0=x^1=0$ we have $x^0=(m+1)/2,y_0=(m-1)/2$ so
that $m\in2\NN+1$.
This proves (b).

\subhead 1.6\endsubhead
In this subsection we assume that $W$ is of type $A_{n-1}$ for some
$n\ge1$. Then $\Irr(W)$ is indexed as in \cite{L84} by the
various $X\sub\ZZ_{>0}$ as in the beginning of 1.4. By the
formulas for $D_E(u)$ in \cite{L84,p.358} we see that if $E\in\Irr(W)$
corresponds to $X$ then $\g_E$ is equal to $\a$ in 1.4(a). Then
1.3(i),(ii) follow from 1.4(b).

We see that the condition that $\Irr_{ssp}(W)\ne\em$ is that
$n=(k^2+k)/2$ for some $k\ge1$; for such $n$ the unique element
$E\in\Irr_{ssp}(W)$ corresponds to $X=\{1,3,\do,2k-1\}$; we have
$c_E=1$,
$$P_E(u)=(1+u)^k(1+u^3)^{k-1}(1+u^5)^{k-2}\do(1+u^{2k-1})/(1+u).$$
We see that 1.3(iii) holds. 

\subhead 1.7\endsubhead
In this subsection we assume that $W$ is of type $B_n$ for some
$n\ge2$. Then $\Irr(W)$ is indexed as in \cite{L84} by the various
$(X,Y)$ as in the beginning of 1.5 such that $\sha(X)=\sha(Y)+1$. By the
formulas for $D_E(u)$ in \cite{L84,p.359} we see that if $E\in\Irr(W)$
corresponds to $(X,Y)$ then $\g_E$ is equal to $\a$ in 1.5(a). Then
1.3(i),(ii) follow from 1.5(b).

We see that the condition that
$\Irr_{ssp}(W)\ne\em$ is that $n=k^2+k$ for some $k$; for such $n$
the unique element $E\in\Irr_{ssp}(W)$ corresponds to
$$(X,Y)=(\{0,2,4,\do,2k\},\{1,3,\do,2k-1\}).$$ We have $c_E=2^k$,
$$P_E(u)=(1+u)^{2k}(1+u^2)^{2k-1}\do(1+u^{2k-1})^2(1+u^{2k}).$$
We see that 1.3(iii) holds. We have $deg(P_E)=2k(k+1)(2k+1)/3$.

\subhead 1.8\endsubhead
In this subsection we assume that $W$ is of type $D_n$ for some
$n\ge4$. Then $\Irr(W)$ is indexed as in \cite{L84} by the various
$(X,Y)$ as in the beginning of 1.5 such that $\sha(X)=\sha(Y)$
except that there are two representations corresponding to any
pair of the form $(X,Y)$ with $X=Y$.

By the
formulas for $D_E(u)$ in \cite{L84,p.359} we see that if $E\in\Irr(W)$
corresponds to $(X,Y)$ then $\g_E$ is equal to $\a$ in 1.5(a). Then
1.3(i),(ii) follow from 1.5(b).

We see that the condition that $\Irr_{ssp}(W)\ne\em$ is that $n=k^2$
for some $k$; for such $n$ the unique element $E\in\Irr_{ssp}(W)$
corresponds to
$$(X,Y)=(\{0,2,4,\do,2k-2\},\{1,3,\do,2k-1\}).$$
We have $c_E=2^{k-1}$,
$$P_E(u)=(1+u)^{2k-1}(1+u^2)^{2k-2}\do(1+u^{2k-2})^2(1+u^{2k-1}).$$
We see that 1.5(iii) holds. We have $deg(P_E)=k(4k^2-1)/3$.

\subhead 1.9\endsubhead
In this subsection we assume that $W$ is of type $E_6$. Using the table
in \cite{L84,p.363}, we see that 1.3(i),(ii) hold; $\Irr_{ssp}(W)$
consists of the unique $E$ such that $\dim(E)=80$. We have $c_E=6$,
$a_E=7$,
$$P_E(u)=(1+u)^3(1+u^2)^2(1+u^3)^3(1+u+u^2)^2.$$
Hence 1.3(iii) holds. We have $deg(P_E)=20$.

\subhead 1.10\endsubhead
In this subsection we assume that $W$ is of type $E_7$. Using the table
in \cite{L84,p.364,365} we see that 1.3(i),(ii) hold; $\Irr_{ssp}(W)$
consists of the unique $E$ such that $\dim(E)=512$, $a_E=11$. We have
$c_E=2$,
$$P_E(u)=(1+u)^2(1+u^3)^2(1+u^5)(1+u^7)(1+u^9).$$
Hence 1.3(iii) holds. We have $deg(P_E)=29$.

\subhead 1.11\endsubhead
In this subsection we assume that $W$ is of type $E_8$.
Using the table in \cite{L84,p.366-369} we see that 1.3(i),(ii) hold;
$\Irr_{ssp}(W)$ consists of the unique $E$ such that $\dim(E)=4480$.
We have  $c_E=120$, $a_E=16$,
$$P_E(u)=(1+u)^4(1+u^2)^4(1+u^3)^4(1+u+u^2)^4(1+u+u^2+u^3+u^4)^2.$$
Hence 1.3(iii) holds. We have $deg(P_E)=40$.

\subhead 1.12\endsubhead
In this subsection we assume that $W$ is of type $F_4$.
Using the table in \cite{L84,p.371} we see that 1.3(i),(ii) hold;
$\Irr_{ssp}(W)$ consists of the unique $E$ such that $\dim(E)=12$.
We have  $c_E=24$, $a_E=4$,
$$P_E(u)=(1+u)^4(1+u^2)^2(1+u+u^2)^2.$$
Hence 1.3(iii) holds.  We have $deg(P_E)=12$.

\subhead 1.13\endsubhead
In this subsection we assume that $W$ is of type $G_2$.
Using the table in \cite{L84,p.372} we see that 1.3(i),(ii) hold;
$\Irr_{ssp}(W)$ consists of the unique $E$ such that $\dim(E)=2$,
$c_E=6$. We have $a_E=1$,
$$P_E(u)=(1+u)^2(1+u+u^2).$$
Hence 1.3(iii) holds. We have $deg(P_E)=4$.

This completes the proof of Theorem 1.3.

\subhead 1.14\endsubhead
In the case where $W$ is irreducible of type $\ne A,E_7$,
the condition
that $W$ is superspecial is equivalent to the following condition:

(a) There exists $E\in\Irr_{sp}(W)$ such that for any
$E'\in\Irr_{sp}(W)-\{E\}$ we have $c_{E'}<c_E$.
\nl
(Then $E$ is unique and is equal to $E_W$.)

\subhead 1.15\endsubhead
It may happen that some non-special $E\in\Irr(W)$ satisfy $\g_E=r$.
For $W$ of type $E_8$, the representations $7168_w,2688_y$ (notation
of \cite{L84}) are such examples; for $W$ of type $F_4$, the
representations $4_1,16_1$ (notation of \cite{L84}) are such examples.
But for other irreducible $W$ there are no such examples.

\subhead 1.16\endsubhead
We can write $W=W_1\T W_2\T\do\T W_e$
where $W_1,W_2,\do,W_e$ are irreducible Weyl groups. If
$E_i,E'_i$ are in $\Irr(W_i)$ ($i=1,\do,e$) then
$E_1\bxt E_2\bxt\do\bxt E_e$, $E'_1\bxt E'_2\bxt\do\bxt E'_e$
are in the same family of $\Irr(W)$ if and only if $E_i,E'_i$ are in
the same family of $\Irr(W_i)$ for $i=1,\do,e$;
we have $E_1\bxt E_2\bxt\do\bxt E_e\in \Irr_{sp}(W)$
if and only if
$E_i\in \Irr_{sp}(W_i)$ for $i=1,\do,e$;
we have $E_1\bxt E_2\bxt\do\bxt E_e\in \Irr_{ssp}(W)$ if and only if
$E_i\in \Irr_{ssp}(W_i)$ for $i=1,\do,e$.

\subhead 1.17\endsubhead
Here is the list of superspecial Weyl groups $W$ that are irreducible
or $\{1\}$.

$A_{(k^2+k)/2-1}$ ($k\in\{1,2,3,\do\}$); 

$B_{k^2+k}$ ($k\in\{1,2,3,\do\}$); 

$D_{k^2}$, $k\in\{2,3,4,\do\}$; 

$E_6,E_7,E_8,F_4,G_2$.

We note that in each case (except in type $A$ and $E_7$) we have that $r$ is
even.

\head 2. The representations $Z_W,Z'_W$\endhead
\subhead 2.1\endsubhead
For any $I'\sub I$ we denote by $W_{I'}$
the subgroup of $W$ generated by $I'$; this is again a Weyl group.
For $E'\in\Irr(W_I)$ we recall that the truncated induction
$J_{W_I}^W(E')$ is the representation of $W$ in which the
the multiplicity of any $E\in\Irr(W)$ is equal to the multiplicity
of $E$ in the ordinary induction $\ind_{W_I}^W(E')$ if $a_E=a_{E'}$
and is $0$ if $a_E\ne a_{E'}$.

\subhead 2.2\endsubhead
In this subsection we assume that $W$ is irreducible or
$\{1\}$, superspecial, simply laced.
We will associate to $W$ a representation $Z_W$ of $W$ of the
form $Z_W=E_1\op E_2\op\do\op E_t$ where $E_1,E_2,\do,E_t$ are
distinct irreducible representations in $\cf_W$ satisfying
the identity
$$c_{E_1}\i+c_{E_2}\i+\do+c_{E_t}\i=1\tag a$$
and (in the case where $W$ is of type $\ne A$) the identity
$$(-1)^{b_{E_1}}c_{E_1}\i+(-1)^{b_{E_2}}c_{E_2}\i+\do+
(-1)^{b_{E_t}}c_{E_t}\i=0\tag a1$$
where $b_E\in\NN$ is defined as in \cite{L84,(4.1.2)}.

If $W$ is of type $A$ (this includes the case $W=\{1\}$)
there is a unique choice for such $Z_W$
namely $Z_W=E_W$ (recall that $c_{E_W}=1$).

If $W$ is of type $E_7$ there is again a unique choice for
such $Z_W$
namely $Z_W=E_1\op E_2$ where $E_1,E_2$ are the two objects of $\cf_W$.
The identities (a),(a1) are now $2\i+2\i=1$. $2\i-2\i=0$.

We now assume that $W$ is of type $\ne A,E_7$.

(b) There is a unique choice of $i\in I$ such that
$I-\{i\}=I'\sqc I''$, $W_{I-\{i\}}=W_{I'}\T W_{I''}$ 
with $W_{I'}$ superspecial, irreducible or $\{1\}$ and with $W_{I''}$
a product of Weyl groups of type $A$ and
such that $E_W$ appears with nonzero multiplicity in
$$J_{W_{I-\{i\}}}^W(E_{W_{I'}}\bxt\sg_{W_{I''}}).$$

We define
$$Z_W=J_{W_{I-\{i\}}}^W(Z_{W_{I'}}\bxt\sg_{W_{I''}}).\tag c$$
(We can assume that $Z_{W_{I'}}$ is known by induction.)
We now describe $Z_W$ in the various cases.

If $W$ is of type $D_{k^2}$ (with $k\ge2$) we have that $W_{I'}$ is of
type $D_{(k-1)^2}$ (if $k\ge3$), $I'=\em$ (if $k=2$) and $W_{I''}$ is of
type $A_{2k-2}$ (if $k\ge3$) and of type $A_1\T A_1\T A_1$ (if $k=2$).

We have $Z_W=\op_\s E_\s$ where $\s$ runs over all permutations of
$1,2,3,\do,2k-2$ which preserve each of the unordered pairs
$(1,2),(3,4),\do,(2k-3,2k-2)$ and for such $\s$, $E_\s\in\Irr(W)$
 corresponds as in 1.8 to 
$$(X,Y)=(\{0,\s(2),\s(4),\do,\s(2k-2)\},\{\s(1),\s(3),\do,\s(2k-3),2k-1\}).$$
Note that $c_{E_\s}=2^{k-1}$ for any $\s$ hence $\sum_\s c_{E_\s}\i=1$.
The identity (a) is now $2^{-k+1}+2^{-k+1}+\do+2^{-k+1}=1$
(the sum has $2^{k-1}$ terms.)

For $E_\s,(X,Y)$ as above we have
$(-1)^{b_{E_\s}}=(-1)^{\sum_{j\in Y}j}h$
where $h=\pm1$ is independent of $\s$ (see \cite{L79a});
hence the identity (a1) holds.

If $W$ is of type $E_6$ we have $I'=\em$ and $W_{I''}$ of type
$A_2\T A_2\T A_1$. We have $Z_W=80_7+60_8+10_9$
(notation of \cite{L15,4.4}. We have
$c_{80_7}=6,c_{60_8}=2,c_{10_9}=3$.
The identities (a),(a1) are now $6\i+2\i+3\i=1$,
$6\i-2\i+3\i=0$.

If $W$ is of type $E_8$ we have $I'=\em$ and $W_{I''}$ of type
$A_4\T A_3$. We have
$$Z_W=4480_{16}+3150_{18}+4200_{18}+420_{20}+7168_{17}+
1344_{19}+2016_{19}$$
(notation of \cite{L15,4.4}). We have

$c_{4480_{16}}=120,c_{3150_{18}}=6,c_{4200_{18}}=8,
c_{420_{20}}=5$,

$c_{7168_{17}}=12,c_{1344_{19}}=4,c_{2016_{19}}=6$.
\nl
The identities (a),(a1) are
now

$120\i+6\i+8\i+5\i+12\i+4\i+6\i=1$,

$120\i+6\i+8\i+5\i-12\i-4\i-6\i=0$.

In each case, $Z_W$ is a constructible representation of $W$
(or a cell in the sense of \cite{L82}); moreover each irreducible
component of $Z_W$ is $2$-special (in the sense if \cite{L15}).

\subhead 2.3\endsubhead
In this subsection we assume that $W$ is irreducible,
superspecial, not simply laced.
We will associate to $W$ two representations $Z_W,Z'_W$ of $W$ of the
form

$Z_W=E_1\op E_2\op\do\op E_t$, $Z'_W=E'_1\op E'_2\op\do\op E'_t$
\nl
where $E_1,E_2,\do,E_t$ (resp. $E'_1,E'_2,\do,E'_t$) are
distinct irreducible representations in $\cf_W$ satisfying
$$c_{E_1}\i+c_{E_2}\i+\do+c_{E_t}\i=1\tag a$$
$$(-1)^{b_{E_1}}c_{E_1}\i+(-1)^{b_{E_2}}c_{E_2}\i+\do
+(-1)^{b_{E_t}}c_{E_t}\i=1\tag a1$$ 
and $c_{E_1}=c_{E'_1},\do,c_{E_t}=c_{E'_t}$,
$b_{E_1}=b_{E'_1},\do,b_{E_t}=b_{E'_t}$.

Now statement 2.2(b) remains true in our case except when
$W$ is of type $B_2,G_2$ or $F_4$ in which case 2.2(b) is true if
one replaces ``unique choice of $i\in I$'' by ``exactly two choices
of $i\in I$''.

If $W$ is of type $B_2,G_2$ or $F_4$ we define $Z_W$ as in 2.2(c)
using one of the two choices of $i$ as above; in these cases we have
$I'=\em$. The same definition using the second choice of $i$ gives
a second representation denoted by $Z'_W$.
If $W$ is of type $B_2$, we have $Z_W=2_1+1_2$,
$Z'_W=2_1+1_2$ (notation of \cite{L15, 4.4}; the $1_2$
in $Z_W$ is different from the $1_2$ in $Z'_W$. 
Then (a),(a1) become  $2\i+2\i=1$, $2\i-2\i=0$.

If $W$ is of type $G_2$, we have $Z_W=2_1+2_2+1_3$
$Z'_W=2_1+2_2+1_3$ (notation of \cite{L15,4.4}); the $1_3$
in $Z_W$ is different from the $1_3$ in $Z'_W$.  
Then (a),(a1)  become $6\i+2\i+3\i=1$, $6\i-2\i+3\i=0$.

If $W$ is of type $F_4$, we have

$Z_W=12_4\op6_6\op9_6\op4_7\op16_5$,

$Z'_W=12_4\op6_6\op9_6\op4_7\op16_5$, (notation
of \cite{L15,4.4}; the $9_6,4_7$ in
$Z_W$ are different from the $9_6,4_7$ in $Z'_W$.
Then (a),(a1) become
$24\i+3\i+8\i+4\i+4\i=1$, $24\i+3\i+8\i-4\i-4\i=0$.

When $W$ is of type $B_{k^2+k},k\ge2$, statement 2.2(b)
remains true.
We define
$$Z_W=J_{W_{I-\{i\}}}^W(Z_{W_{I'}}\bxt\sg_{W_{I''}}),$$
$$Z'_W=J_{W_{I-\{i\}}}^W(Z'_{W_{I'}}\bxt\sg_{W_{I''}}),$$
(We can
assume that $Z_{W_{I'}},Z'_{W_{I'}}$ are known by induction.)

We have $Z_W=\op_\s E_\s$ where $\s$ runs over all permutations of
$1,2,3,\do,2k$ which preserve each of the unordered pairs
$(1,2),(3,4),\do,(2k-1,2k)$ and for such $\s$,
$E_\s\in\Irr(W)$
corresponds as in 1.7 to 
$$(X,Y)=(\{0,\s(2),\s(4),\do,\s(2k)\},\{\s(1),\s(3),\do,\s(2k-1)\}).$$
Note that $c_{E_\s}=2^k$ for any $\s$ hence $\sum_\s c_{E_\s}\i=1$.
For $E_\s,(X,Y)$ as above we have
$(-1)^{b_{E_\s}}=(-1)^{\sum_{j\in Y}j}h$
where $h=\pm1$ is independent of $\s$ (see \cite{L79a});
hence the identity (a1) holds.

We have $Z'_W=\op_\s E'_\s$ where $\s$ runs over all permutations of
$0,1,2,3,\do,2k-1$ which preserve each of the unordered pairs
$(0,1),(2,3),\do,(2k-2,2k-1)$ and for such $\s$, $E'_\s\in\Irr(W)$
corresponds as in 1.7 to 
$$(X,Y)=(\{\s(0),\s(2),\s(4),\do,\s(2k-2),2k\},\{\s(1),\s(3),\do,\s(2k-1)\}).$$

In each case, $Z_W,Z'_W$ are constructible representations of $W$
(or cells in the sense of \cite{L82}) such that
$Z'_W=Z_W\ot\sg_W$; their irreducible
components are $2$-special (in the sense if \cite{L15}).

\head 3. The noncrystallographic case\endhead
\subhead 3.1 \endsubhead
In this section we assume that $W$ is an irreducible
noncrystallographic finite Coxeter group with set $I$ of simple
reflections. Now $\Irr(W)$ is defined as in 0.1;
the generic degree $D_E(u)$ for
$E\in\Irr(E)$ is defined as in 1.1. (We now have
$D_E(u)\in\RR[u]$, see \cite{AL,p.202}, \cite{L82}.)
From the
formulas for $D_E(u)$ in type $H_4$ in \cite{AL} one can verify
that an equality like 1.1(a) still
holds except that now $c_E$ is only an algebraic integer in
$\RR_{>0}$ and $P_E(u)$
is now only in $\RR[u]$ with leading coefficient $1$; a similar
property holds in the cases $\ne H_4$. Then $a_E\in\NN$ is defined
as in 1.1. The definition of families in \cite{L79} and that
of $\Irr_{sp}(W)$ extend in an obvious way to our case.
(The representations in $\Irr_{sp}(W)$ are described explicitly
in type $H_4$ in \cite{AL,\S5}.)
For $E\in\Irr(W)$ we define $\g_E\in\NN$ as in 1.2.

\proclaim{Theorem 3.2} (i) For any $E\in\Irr(W)$ we have $\g_E\le\sha(I)$.

(ii) The set $\Irr_{ssp}(W):=\{E\in\Irr_{sp}(W);\g_E=\sha(I)\}$ consists
of exactly one element.

(iii) If $E\in\Irr_{ssp}(W)$ then $P_E(u)\in\RR_{\ge0}[u]$. It has
degree $2a_E+\sha(I)$.
\endproclaim
This can be verified using the known results on $D_E(u)$.
Assume first that $W$ is of type $H_4$. Now $\Irr_{ssp}(W)$ consists of
the unique $E$ such that $\dim(E)=24$, $a_E=6$. We have
$c_E=120/(13-8\l)=120(5+8\l)$,
$$P_E(u)=(u+1)^4(u^2+u+1)^2(u^2+\l u+1)^2(u^2-\ti\l u+1)^2$$
where $$\l=(1+\sqrt5)/2\in\RR_{>0},-\ti\l=(\sqrt5-1)/2\in\RR_{>0}.$$

Assume next that $W$ is of type $H_3$. Now
$\Irr_{ssp}(W)$ consists of the unique $E$ such that $\dim(E)=4$,
$a_E=3$. We have $c_E=2$, $P_E(u)=(1+u)(1+u^3)(1+u^5)$.

If $W$ is a dihedral group of order $2p$, $p=5$ or $p\ge7$ then
$\Irr_{ssp}(W)$ consists of the unique $E$ such that $\dim(E)=2$,
$a_E=1$. We have $$c_E=p/((1-\x)(1-\x\i))=\prod_{t=2}^{p-2}(1-\x^t),$$
$P_E=(1+u)^2(1+(\x+\x\i)u+u^2)$ where $\x=e^{2\pi\sqrt{-1}/p}$.

\head 4. Superspecial conjugacy classes in $W$\endhead
\subhead 4.1 \endsubhead
In this section we assume that $W$ is an
irreducible superspecial Weyl group. Let $E=E_W$.

Until the end of 4.3 we assume also that
the opposition $op:I@>>>I$ is the identity map.

Let $G,G(F_q),\r_E,\r'$ be as in 0.1.
Recall that $\r'$ is a unipotent cuspidal representation of
$G(F_q)$, say over $\bbq$.
For any $w\in W$ let $X_w$ be the subvariety of the flag manifold
of $G$ defined in \cite{DL} (it consists of Borel subgroups
which are in
relative position $w$ with their transform under the Frobenius
map. Let $H^i_c(X_w)$ be the $i$-th $l$-adic
cohomology with compact support of $X_w$ viewed as a representation
of $G(F_q)$. Let $c(W)$ be the set of all $w\in W$ such that $\r'$
appears with nonzero multiplicity in the virtual representation
$\sum_iH^i_c(X_w)$ of $G(F_q)$.
Let $M(W)$ be the minimum of the legths of various elements
in $c(W)$ and let $c_{min}(W)$ be the set of elements of $c(W)$
of length $M(W)$.
The following result can be deduced from \cite{L02,2.18}.

\proclaim{Theorem 4.2} There is a unique conjugacy class $C_W$
in $W$ such that $c_{min}(W)\sub C_W$.
\endproclaim
The conjugacy class $C_W$ is said be the superspecial conjugacy
class of $W$. It is an elliptic conjugacy class. 

\subhead 4.3\endsubhead
We describe $C_W$ in each case (we use \cite{L02}).
We also describe in each case the number $M(W)$ (we use \cite{GP}).

If $W$ is of type $B_{k^2+k},k\ge1$ viewed in an obvious way as a
subgroup of $S_{2(k^2+k)}$ then $C_W$ is the elliptic conjugacy
class with cycle type $4+8+12+\do+4k$. We have $M(W)=k(k+1)(2k+1)/3$.

If $W$ is of type $D_{k^2}$, $k\ge2$ even, viewed in an obvious way as a
subgroup of $S_{2k^2}$ then $C_W$ is the elliptic conjugacy
class with cycle type $2+6+10+\do+(4k-2)$. We have $M(W)=2k(k^2-1)/3$.

If $W$ is of type $E_8$ then $C_W$ consists of elements with
characteristic polynomial $(u^2-u+1)^4$ in the reflection
representation. We have $M(W)=deg(P_E(u))=40$.
We have $c(W)=c_{min}(W)=C_W$, $\sha(C_W)=\dim(E)$.

If $W$ is of type $E_7$ then $C_W$ is the Coxeter conjugacy class.
We have $M(W)=7$.

If $W$ is of type $F_4$ then $C_W$ consists of elements
with characteristic polynomial $(u^2+1)^2$ in the reflection
representation. We have $M(W)=deg(P_E(u))=12$.
We have $c(W)=c_{min}(W)=C_W$, $\sha(C_W)=\dim(E)$.

If $W$ is of type $G_2$ then $C_W$ consists of elements
with characteristic polynomial $u^2+u+1$ in the reflection
representation. We have $M(W)=deg(P_E(u))=4$.
We have $c(W)=c_{min}(W)=C_W$, $\sha(C_W)=\dim(E)$.

\subhead 4.4\endsubhead
We now assume that the opposition $op:I@>>>I$ is not the
identity map. It induces an involution $op:W@>>>W$.
Then results similar to 4.2 hold (we use \cite{L02,2.19}).
They associate to $W$ a
twisted conjugacy class in $W$ (an  orbit of the $W$-action
$x:w\m xw op(x)\i$ on $W$).


\widestnumber\key{ABCD}
\Refs
\ref\key{AL}\by D.Alvis and G.Lusztig \paper  The representations and generic degrees of the Hecke algebras of type $H_4$, J. reine und angew.math.\vol336\yr1982\pages201-212\moreref Erratum\vol449\yr1994\pages217-281\endref
\ref\key{DL}\by P.Deligne and G.Lusztig\paper Representations of
reductive groups over finite fields\jour Ann. Math.\vol103\yr1976\pages
103-161\endref
\ref\key{GP}\by M.Geck and G.Pfeiffer\book Characters of finite
Coxeter groups and Iwahori-Hecke algebras\publ Clarendon Press,
Oxford\yr2000\endref
\ref\key{L77}\by G.Lusztig\paper Irreducible representations of finite
classical groups\jour Inv. Math.\vol43\yr1977\pages125-175\endref
\ref\key{L79}\by G.Lusztig\paper Unipotent representations of a finite
Chevalley group of type $E_8$\jour Quart. J. Math.\vol30\yr1979\pages
315-338\endref
\ref\key{L79a}\by G.Lusztig\paper A class of irreducible
representations of a Weyl group\jour Proc. Kon. Nederl. Akad.(A)\vol82
\yr1979\pages 323-335\endref
\ref\key{L82}\by G.Lusztig\paper A class of irreducible representations
of a Weyl group II\jour Proc. Kon. Nederl. Akad.(A)\vol85\yr1982\pages
219-226\endref
\ref\key{L84}\by G.Lusztig\book Characters of reductive groups over a
finite field\bookinfo Ann.Math.Studies 107\publ Princeton U.Press\yr
1984\endref
\ref\key{L02}\by G.Lusztig\paper Rationality properties of unipotent
representations\jour J. Alg.\vol258\yr2002\pages1-22\endref
\ref\key{L15}\by G.Lusztig\paper On conjugacy classes in a
reductive group\inbook Representations of Reductive 
Groups\bookinfo Progr.in Math. 312\publ Birkh\"auser\yr2015\pages
333-363\endref
\endRefs
\enddocument